\documentclass[12pt]{article}

\usepackage{latexsym,amssymb,amscd,epic,eepic,epsfig,graphics}
\usepackage{amsbsy,amsmath,amsfonts}
\usepackage{amsthm}
\usepackage[francais,english]{babel}
\usepackage[all]{xy}

\newtheorem{theor}{Th\'eor\`eme}
\newtheorem{thm}{Th\'eor\`eme}
\newtheorem{lemme}{Lemme}
\newtheorem{cor}{Corollaire}
\newtheorem{prop}{Proposition}
\newtheorem{defi}{D\'efinition}
\newtheorem{rem}{Remarque}
\newtheorem{notation}{Notations}
\newtheorem*{algo}{Algorithme}
\newtheorem*{Merci}{Remerciements}

\newenvironment{demo}
{\begin{proof}[D\'emonstration]} {\end{proof}}

\newenvironment{sdemo}[1]
{\begin{proof}[#1]} {\end{proof}}

\newcommand{\ov}[1]{\overline{#1}}

\newcommand{\N}{\mathbb N}
\newcommand{\Z}{\mathbb Z}
\newcommand{\R}{\mathbb R}

\newcommand{\Lm}{\mathbb{L}_m}
\newcommand{\Gm}{\mathcal{G}_m}
\newcommand{\GZm}{\mathcal{G}(\mathbb{Z}^m)}
\newcommand{\Pm}{\mathcal{P}(\mathbb{L}_m)}

\newcommand{\BZn}{B_{\mathcal{A}}(n)}
\newcommand{\B}[1]{B_{S}\left(#1 \right)}

\newcommand{\cln}[1]{\mathcal{N} \left( #1 \right) }
\newcommand{\dimi}{\underline{\dim}_M}
\newcommand{\dims}{\overline{\dim}_M}

\begin{document}

\selectlanguage{francais}


\date{20 f\'evrier 2005}


\author{Luc GUYOT \footnote{Ce travail est soutenu par le Fond National Suisse, No.~PP002-68627.}
\footnote{Section de Math\'ematiques de Gen\`eve \newline  2-4, rue du Li\`evre \newline Case Postale 240 \newline
CH-1211 Gen\`eve 24 \newline \newline E-mail : Luc.Guyot@math.unige.ch}}


\title{ Estimations de dimensions de Minkowski dans l'espace des groupes marqu\'es }


\maketitle

\begin{abstract}
Dans cet article, on montre que l'espace des groupes marqu\'es est un sous-espace
 ferm\'e d'un ensemble de Cantor dont la dimension de Hausdorff est infinie.
 On prouve que la dimension de Minkowski de cet espace est infinie en exhibant des sous-ensembles de groupes
 marqu\'es \`a petite simplification dont les dimensions de Minkowski sont arbitrairement grandes.
 On donne une estimation des dimensions de Minkowski de sous-espaces de groupes \`a un relateur.
 On d\'emontre enfin que les dimensions de Minkowski du sous-espace des groupes
 commutatifs marqu\'es et d'un ensemble de Cantor
  d\'efini par Grigorchuk sont nulles.
\end{abstract}

\selectlanguage{english}

\begin{abstract}
In this article we show that the space of marked groups is a closed
subspace of a Cantor space with infinite Hausdorff dimension. We
prove that the Minkowski dimension of this space is infinite by
exhibiting subsets of marked groups with small cancellation the
dimensions of which are arbitrarly large. We give estimates of the
Minkowski dimensions of
 of subsets of marked groups with one relator. Eventually, we prove that the Minkowski dimensions of the subspace of
abelian marked groups and a Cantor space defined by Grigorchuk are zero.
\end{abstract}

\selectlanguage{francais}

\section{Introduction} \label{intro}%

L'\'etude de la g\'en\'ericit\'e de diff\'erentes classes de groupes a donn\'e lieu \`a de nombreux travaux
depuis le th\'eor\`eme de g\'en\'ericit\'e des groupes hyperboliques \'enonc\'e par Gromov \cite{grom2}.
Un nouvel aspect dans la caract\'erisation de cette g\'en\'ericit\'e a \'et\'e d\'evelopp\'e par Champetier \cite{champ1}
en consid\'erant l'espace topologique des groupes marqu\'es \`a $m$ g\'en\'erateurs et les cat\'egories de Baire de parties
sp\'ecifiques de cet espace.

Dans l'id\'ee de caract\'eriser cette g\'en\'ericit\'e d'un point de vue m\'etrique et
de mesurer l'importance relative de certaines classes de groupes, on s'int\'eresse
 ici \`a l'estimation des dimensions de Minkowski et de Hausdorff de l'espace des groupes marqu\'es \`a $m$ g\'en\'erateurs
  et de certaines de ses parties, qu'on munit de la m\'etrique employ\'ee par Champetier \cite{champ1}.

On d\'ecrit dans un premier temps un plongement isom\'etrique naturel de l'espace des groupes marqu\'es dans un Cantor dont on montre
que la dimension de Hausdorff est infinie.
On consid\`ere ensuite la partie $PS=PS(m,k,\lambda)$ form\'ee des groupes \`a $k$ relateurs de m\^eme longueur v\'erifiant la
condition de petite simplification $C'(\lambda)$~:
 \begin{theor}
 Lorsque $m \ge 2$ et $\lambda \in \rbrack 0,\frac{1}{6} \rbrack$, on a l'encadrement des dimensions inf\'erieure et
 sup\'erieure de Minkowski suivant

$$
k \log_2(2m-1) \le \dimi PS \le \dims PS \le \frac{k}{1-3\lambda}\log_2(2m-1).
$$
 \end{theor}

\begin{cor}
La dimension de Minkowski inf\'erieure de l'espace des groupes marqu\'es \`a $m$ g\'en\'erateurs ($m \ge 2$) est infinie.
\end{cor}

On introduit apr\`es cela la partie $UR=UR(m,q)$ des groupes \`a un relateur et dont le relateur est
une puissance $q$-\`eme, partie pour laquelle on montre le
\begin{theor} Lorsque $q \ge 2$, on a l'encadrement
$\frac{\log_2(2m-1)}{q} \le \dimi UR \le \dims UR \le \frac{\log_2(2m-1)}{q-1}.$
\end{theor}
Dans le cas de $m=4$ g\'en\'erateurs, on s'int\'eresse au sous-espace $\mathfrak{B}$ d\'efini par Grigorchuk \cite{gri2}
 dont on rappelle en d\'etail la construction au chapitre \ref{secdimB} et pour lequel on montre
\begin{theor}
La dimension de Minkowski sup\'erieure de l'espace $\mathfrak{B}$ est nulle.
\end{theor}
On montre enfin le
\begin{theor}
La dimension de Minkowski sup\'erieure du sous-espace des groupes commutatifs est nulle.
\end{theor}

 Les chapitres \ref{secespace} et \ref{secdim} sont deux chapitres pr\'eliminaires o\`u toutes les d\'efinitions utiles
 regardant l'espace des groupes marqu\'es et les dimensions m\'etriques de Minkowski et de Hausdorff ont \'et\'e
 rassembl\'ees ainsi que quelques exemples. Dans le chapitre \ref{secPLm} on prouve, lorsque $m \ge 2$,
  que l'espace $\mathcal{P}(\Lm)$ des parties de $\Lm$ dans lequel se plonge
   isom\'etriquement $\mathcal{G}_m$ a une dimension de Hausdorff infinie. On y montre
 \'egalement que la dimension sup\'erieure de Minkowski de $\mathcal{G}_1$ est nulle. Les chapitres \ref{secdimPS} \`a
 \ref{secdimGZ} pr\'esentent dans l'ordre les d\'emonstrations des quatre r\'esultats principaux \'enonc\'es
 dans cette introduction.

\begin{Merci}
 Je tiens \`a remercier tr\`es chaleureusement Roland Bacher, autant pour ses id\'ees pr\'ecieuses que ses remarques
 judicieuses. Je remercie \'egalement Goulnara Arjantseva et Pierre De La Harpe pour leurs patientes et minutieuses relectures, ainsi que
  Fr\'ed\'eric Mouton pour ses tr\`es opportunes suggestions.
\end{Merci}

\section{Les espaces $\mathcal{G}(G)$ et $\mathcal{P}(G)$} \label{secespace} %
Soit $\Lm$ le groupe libre de base $S=\{e_1,\ldots,e_m\}.$ A chaque sous-groupe distingu\'e de $\Lm$
 correspond un quotient marqu\'e de $\Lm$, c'est-\`a-dire un groupe
muni d'un syst\`eme de g\'en\'erateurs privil\'egi\'e qui est l'image du syst\`eme $S$
 par l'application quotient. L'ensemble $\mathcal{G}_m$ des sous-groupes distingu\'es de $\Lm$, lorsqu'il est
 muni d'une topologie m\'etrisable dite topologie de Cayley, est un espace compact appel\'e espace des groupes
 marqu\'es \`a $m$ g\'en\'erateurs.

La topologie de Cayley sur $\mathcal{G}_m$,
est induite par une topologie m\'etrisable sur l'ensemble $\mathcal{P}(\Lm)$ des parties de $\Lm$ qui fait
de $\mathcal{P}(\Lm)$ un espace de Cantor (un espace de Cantor est un espace topologique compact,
 totalement discontinu et sans points isol\'es;
un tel espace est hom\'eomorphe \`a l'ensemble triadique de Cantor).
On donne ci-dessous une construction plus g\'en\'erale qui permet
de d\'efinir une topologie compacte et m\'etrisable sur l'ensemble des parties d'un groupe de type fini ainsi que sur
l'ensemble de ses sous-groupes distingu\'es.

Soit $G$ un groupe de type fini muni d'un syst\`eme ordonn\'e de g\'en\'erateurs  $X=(g_1,\ldots,g_m)$.
La longueur $\vert g \vert_{X}$ d'un \'el\'ement $g \in G$ relativement \`a $X$
 est la longueur d'un mot irr\'eductible le plus court en les lettres $X
\cup X^{-1}$ qui repr\'esente $g$. On d\'esigne par $B_{X}(r)$ l'ensemble des \'el\'ements de
$G$ de longueur inf\'erieure ou \'egale \`a $r$.

 Si $R$ est une partie de $G$, on d\'esigne par $\cln{R}$ le plus petit sous-groupe distingu\'e de
$G$ contenant $R$ appel\'e cl\^oture normale de $R$. On \'ecrit $<R>$ pour le sous-groupe de $G$ engendr\'e par $R$. Si $A$ est un ensemble fini,
 on d\'esigne par $\vert A \vert$ le cardinal de $A$.

\begin{defi}
Une distance $d$ sur un ensemble $E$ v\'erifie l'in\'egalit\'e ultra-m\'etrique si
$$
d(x,z) \le \max \left\{ d(x,y),d(y,z) \right\},
\mbox{ pour tout $x,y$ et $z$ de $E$}.
$$
\end{defi}

Soit  $\mathcal{P}(G)$ l'ensemble des parties de $G$.
 On d\'efinit sur $\mathcal{P}(G)$ la distance ultra-m\'etrique $d_{\mathcal{P}}$
  \`a partir de la ``valuation''
$$\nu(A,B)= \max \left\{ n \in \N \cup \{\infty\} \,\vert \, B_{X}(n) \cap A= B_{X}
(n) \cap B \right\}$$
 pour des parties $A$ et $B$ de $G$. On pose alors
$d_{\mathcal{P}}(A,B)=2^{-\nu(A,B)}.$
Il est imm\'ediat de v\'erifier que la topologie induite par $d_{\mathcal{P}}$ sur $\mathcal{P}(G)$
 est aussi celle de la topologie produit sur $\{0,1\}^G$ o\`u $\{0,1\}$ est muni de la topologie discr\`ete.
On montre sans peine que
 $\left( \mathcal{P}(G), d_{\mathcal{P}} \right)$ est un espace de Cantor si $G$ est infini et un espace
 discret fini sinon.
\begin{defi}
On d\'esigne par $\mathcal{G}(G)$ l'ensemble des sous-groupes distingu\'es de $G$.
On consid\`ere sur $\mathcal{G}(G)$ la m\'etrique $d_{\mathcal{G}}$ induite par
$d_{\mathcal{P}}$. L'espace m\'etrique $\mathcal{G}(G)$ est appel\'e espace des quotients marqu\'es de $G$.
Lorsque $G=\Lm$, on \'ecrit $\mathcal{G}_m=\mathcal{G}(\Lm)$.
\end{defi}

L'espace $\mathcal{G}(G)$ est une partie ferm\'ee de $\mathcal{P}(G)$ qui poss\`ede des points isol\'es.
Ainsi $\mathcal{G}(G)$ est un espace m\'etrique compact totalement discontinu et donc
 de dimension topologique nulle, voir par exemple ~\cite{champ1,champ2} et~\cite{gri}. Lorsque $G=\Lm$, on \'ecrit
$d_{\mathcal{G}}=d_m$.

\section{Dimensions de Minkowski et de Hausdorff}    \label{secdim}      %

 Les dimensions de Minkowski et de Hausdorff sont des dimensions m\'etriques qui
renseignent sur les possibilit\'es de plongement dans un espace m\'etrique standard tel qu'un espace
euclidien ou hyperbolique. Pour les notions de base concernant
 les dimensions m\'etriques on se r\'ef\`ere au livre de
 Falconer \cite{falc}.

 On d\'esigne par $(E,d)$ un espace m\'etrique.
 Pour toute partie $A$ de $E$, $diam(A)$ est le diam\`etre de $A$.
Un espace pr\'ecompact est un espace m\'etrique qui poss\`ede pour tout $\varepsilon >0$ un recouvrement fini
 par des boules ferm\'ees de rayon $\varepsilon$.

\begin{notation} \label{defirecmin}
Soit $\left( E,d \right)$ un espace m\'etrique pr\'ecompact.
On d\'esigne par $N(E,\varepsilon)$ le minimum des cardinaux des recouvrements
de $E$ par des boules ferm\'ees de rayon $\varepsilon$. On d\'esigne par $P(E,\varepsilon)$
 le nombre maximum de boules ferm\'ees
de rayon $\varepsilon$ disjointes.
\end{notation}

\begin{defi}
Les dimensions de Minkowski inf\'erieure et sup\'erieure d'un espace m\'etrique pr\'ecompact $(E,d)$
 se d\'efinissent respectivement par les formules
$$
\dimi E=\liminf_{\varepsilon \to 0} \frac{\log N(E,\varepsilon)}{\log(1/\varepsilon)},
$$
$$
\dims E=\limsup_{\varepsilon \to 0}\frac{\log N(E,\varepsilon)}{\log(1/\varepsilon)}.
$$
\end{defi}
Ces dimensions m\'etriques sont aussi connues sous les noms de dimensions fractales ou ``box-counting dimensions''
inf\'erieure et sup\'erieure.

\begin{defi}
La dimension de Hausdorff d'un ensemble $A \subset E$ est
$$
\dim_H A =\sup \left\{ s : \mathcal{H}^s(A)>0 \right\}=\sup \left\{
s: \mathcal{H}^s(A)=\infty \right\}$$
$$
= \inf \left\{ s : \mathcal{H}^s(A)< \infty \right\}= \inf \left\{s : \mathcal{H}^s(A)=0 \right\},
$$
o\`u $\mathcal{H}^s$ est la mesure de Hausdorff de dimension $s$
$$
\mathcal{H}^s(A)= \lim_{ \delta \to 0} \inf \left\{ \sum_{i=1}^{\infty}diam(E_i)^s : A \subset \bigcup_{i=1}^{\infty}E_i,
 \, E_i \subset E,\,diam(E_i) \le \delta \right\}.
$$
\end{defi}

Les propri\'et\'es \'el\'ementaires suivantes sont v\'erifi\'ees lorsque $\dim$ d\'esigne soit $\dimi,\dims$ ou $\dim_H$
(voir \cite[Ch.2]{falc}).

Monotonie~:
Si $E_1 \subset E_2$ alors $\dim E_1 \le \dim E_2$. \\
Ensemble fini~:
Si $E$ est fini alors $\dim E=0$.

La dimension de Minkowski sup\'erieure est \it finiment stable \rm, c'est-\`a-dire

$$
\dims  \bigcup_{i=1}^{n}E_i =\underset{1 \le i \le n}{\max} \dims E_i \, , \mbox{ pour } E_i \subset E, i=1,2,\ldots,n,
$$
alors que la dimension de Hausdorff est \it d\'enombrablement stable \rm, c'est-\`a-dire

$$
\dim_H  \bigcup_{i=1}^{\infty}E_i =\underset{i \ge 1}{\sup} \dim_H E_i \, , \mbox{ pour } E_i \subset E, i=1,2,\ldots .
$$
La dimension de Hausdorff d'un ensemble d\'enombrable est donc nulle.
Si $\dim$ d\'esigne l'une des dimensions de Minkowski, alors $\dim A=\dim \overline{A}$ o\`u $\overline{A}$
est l'adh\'erence de $A \subset E.$ Pour cette raison, dimension de Hausdorff et dimensions
de Minkowski peuvent \^etre très diff\'erentes.

\begin{prop} \cite{matt} Pour tout espace m\'etrique pr\'ecompact $E$,
\begin{displaymath} \label{inedim}
\dim_{top} E \le \dim_H E \le \dimi E \le \dims E.
\end{displaymath}
o\`u $\dim_{top}$ d\'esigne la dimension topologique ou dimension de recouvrement de Lebesgue.
\end{prop}

\subsection{Exemples}
Si $E$ est la boule unit\'e d'un espace vectoriel norm\'e de dimension $n$, munie de l'une quelconque de ses normes, alors ses dimensions de Hausdorff
et de Minkowski sont toutes \'egales \`a $n$ qui est aussi sa dimension topologique. Si $E$ est l'espace triadique
de Cantor, construit dans le segment $\lbrack 0,1 \rbrack$ muni de la m\'etrique euclidienne, ses dimensions de Minkowski et de Hausdorff valent toutes trois
$\log(2)/\log(3)$ alors que sa dimension topologique est nulle.\\
Si $E=\lbrack 0,1 \rbrack$, muni de la m\'etrique euclidienne, et si $\alpha$ est un r\'eel strictement positif, alors
\begin{itemize}
\item[$(1)$]
$\dimi A=\dims A=\frac{1}{1+\alpha}$ lorsque $A=\{ \frac{1}{n^{\alpha}} \}_{n \ge 1}$,
\item[$(2)$]
$\dimi A=\dims A=1$ lorsque $A=\{ \frac{1}{\log n } \}_{n \ge 2}$,
\item[$(3)$]
$\dimi A=\dims A=0$ lorsque $A=\{ 2^{-n} \}_{n \ge 0}$.
\end{itemize}
Dans chacun de ces cas, la dimension de Hausdorff est nulle. Ces r\'esultats s'obtiennent directement \`a partir
 de la d\'efinition.
  Le troisi\`eme cas sugg\`ere que la dimension de Minkowski sup\'erieure de l'ensemble des valeurs
 d'une suite  qui converge vers son unique point d'accumulation \`a vitesse exponentielle,  est nulle.
  C'est le cas dans l'espace
 $\mathcal{G}(\Z)$ o\`u $d_1(n\Z,\{0\})=2^{-n}$ et o\`u le sous-groupe distingu\'e trivial
  est l'unique point
 d'accumulation de cet espace. Comme premier exemple de calcul,
  on montre dans la proposition \ref{dimZ} du chapitre suivant
  qu'effectivement $\dims \mathcal{G}(\Z)=0$. De mani\`ere g\'en\'erale, majorer
   la dimension de Minkowski sup\'erieure de l'ensemble des valeurs d'une suite convergeant
   vers son unique point d'accumulation revient estimer la vitesse de convergence de cette suite.

\section{Dimension de Hausdorff de l'espace $\mathcal{P}(\Lm)$} \label{secPLm}%

On consid\`ere, comme dans le chapitre \ref{secespace}, un groupe $G$ muni
d'un syst\`eme  ordonn\'e de g\'en\'erateurs $X$ ainsi que
l'espace m\'etrique $\left( \mathcal{P}(G),d_{\mathcal{P}} \right)$.
On pose $\beta(n)=\vert B_{X}(n) \vert$ et $\sigma(n)=\beta(n+1)-\beta(n).$
La limite $\omega(G,X)=\underset{n \to \infty}{\limsup}\sqrt[n]{\beta(n)}$ est appel\'ee taux de
 croissance exponentiel du groupe $G$ relativement \`a $X$. Si $\omega(G,X)>1$, on dit que
le groupe $G$ est de \it croissance exponentielle.\rm  Cette derni\`ere d\'efinition ne d\'epend pas du choix de
$X$ (voir \cite[chap.VI.C]{DeHa}). On dit que deux points d'un espace m\'etrique sont $\varepsilon$-distinguables
si la distance qui les s\'epare est strictement sup\'erieure \`a $\varepsilon$.

On exprime dans la propositon \ref{propdimMPG} les dimensions de Minkowski de $\mathcal{P}(G)$ en fonction de $\beta$.
Le lemme \ref{lemmecritere} donne une condition suffisante portant sur la croissance de $G$
 pour que l'egalit\'e ait lieu  entre dimension de Hausdorff et dimension de Minkowski inf\'erieure de l'espace
  $\mathcal{P}(G)$. Il en r\'esulte imm\'ediatement que la dimension de Hausdorff de $\mathcal{P}(\Lm)$ est infinie.
La d\'emonstration de ce premier lemme est laiss\'ee au lecteur.

\begin{lemme} \label{lemmerecouvrement}
Soit $\left(E,d \right)$ un espace ultra-m\'etrique compact.
Pour tout $\varepsilon >0$, il
 existe un unique recouvrement fini minimal $\mathcal{F}(E,\varepsilon)$ de $E$
  par des boules ferm\'ees de rayon $\varepsilon$, c'est-\`a-dire tel qu'aucune sous-famille
 propre de $\mathcal{F}(E,\varepsilon)$ n'est un recouvrement de $E$.
 De plus ce recouvrement est une partition et l'on a
 $\vert \mathcal{F}(E,\varepsilon) \vert=N \left( E,\varepsilon \right)=P \left(E, \varepsilon \right)$.
 En outre, $P(E,\varepsilon)$  est \'egal au nombre maximal de points $\varepsilon$-distinguables dans $E$.
\end{lemme}

\begin{lemme} \label{lemmerecouvrement2}
Soit $G$ un groupe de type fini muni d'un syst\`eme ordonn\'e de g\'en\'erateurs $X$.
Si $(E,d)=\left(\mathcal{P}(G),d_{\mathcal{P}} \right)$ ou $\left(\mathcal{G}(G), d_{\mathcal{G}} \right)$ alors,
 pour tout entier $n \ge 0$,
$N(E,2^{-n})$ est le nombre de parties de la boule $B_{X}(n)$ qui s'obtiennent comme
l'intersection d'un \'el\'ement de $E$ avec cette boule.
\end{lemme}

\begin{demo}
Supposons que $ \left( E, d \right) =\left( \mathcal{P}(G), d_{\mathcal{P}} \right)$. On se donne $n \ge 0$.
 Ayant pos\'e $\varepsilon =2^{-n}$, on consid\`ere le sous-ensemble $P_n= \mathcal{P}(B_{X}(n))$ de $E$
form\'e des parties de la boule $B_{X}(n)$.
Alors, l'ensemble des boules centr\'ees en les points de $P_n$ et de rayon $\varepsilon$
est un recouvrement de $E$ par des boules disjointes. C'est donc le recouvrement minimal
$\mathcal{F} \left( E,\varepsilon \right)$ d'apr\`es le lemme \ref{lemmerecouvrement}.
 En effet, deux centres $x_i$ et $x_j$ distincts v\'erifient
$ d(x_i,x_j) > \varepsilon$.
Si les boules associ\'ees n'\'etaient pas disjointes, elles seraient alors \'egales \`a une m\^eme boule
de diam\`etre strictement sup\'erieur \`a $\varepsilon$. Ceci est absurde puisque le diam\`etre d'une boule de rayon $\varepsilon$ d'un espace
ultra-m\'etrique n'exc\`ede pas $\varepsilon$.
Lorsque $(E,d)=\left(\mathcal{G}(G), d_{\mathcal{G}} \right)$, la preuve reste en tout point semblable.
\end{demo}

\begin{prop} \label{propdimMPG}
Les dimensions de Minkowski inf\'erieures et sup\'erieures de l'espace m\'etrique
$\left(\mathcal{P}(G),d_{\mathcal{P}(G)}  \right)$ sont donn\'ees par les formules
$$
\dimi \mathcal{P}(G) =\underset{n \to \infty}{\liminf} \frac{\beta(n)}{n} ,\quad
\dimi \mathcal{P}(G) =\underset{n \to \infty}{\limsup} \frac{\beta(n)}{n}.
$$
\end{prop}

\begin{demo}
Ayant fix\'e $n \ge 0$, on pose $\varepsilon=2^{-n}$.
Par le lemme \ref{lemmerecouvrement2},
$$ N \left(E,\varepsilon \right)=2^{ \vert B_{X}(n) \vert }=2^{\beta(n)}.$$

$$ \mbox{D'o\`u }
\dimi \Pm = \liminf_{ \varepsilon \to 0} \frac{ \log N \left( E, \varepsilon \right)}{ \log(1/ \varepsilon)}=
\liminf_{n \to \infty} \frac{ \beta(n)}{n},
$$
avec une formule analogue pour la dimension de Minkowski sup\'erieure.
\end{demo}

Si $d$ est une distance sur un ensemble $E$, on d\'esigne par $val(d)$ le
 sous-ensemble de $\R^+$ constitu\'e des valeurs prises par $d$ sur $E \times E$.
Pour tout $x \in E$, $\varepsilon>0$, on d\'esigne par $B(x,\varepsilon)$ la boule ferm\'ee de centre $x$ et
 de rayon $\varepsilon$.

Le lemme suivant montre que sous une condition de contr\^ole uniforme des ``$s$-volumes"
$\varepsilon'^{s}N(B(x,\varepsilon),\varepsilon')$ des partitions r\'eguli\`eres des boules,
 dimension de Hausdorff et dimension inf\'erieure de Minkowski d'un espace ultra-m\'etrique compact sont \'egales.
On laisse au lecteur la preuve de ce lemme qui est \'el\'ementaire.

\begin{lemme} \label{lemmecritere}

Soit $\left( E, d \right)$ un espace ultra-m\'etrique compact v\'erifiant les deux hypoth\`eses suivantes

\begin{itemize}

\item[$(1)$]
Pour toute boule $B$ de rayon $r \in val(d)$, $diam(B)=r$,
\item[$(2)$]
Pour tout $s>0$, il existe un nombre r\'eel $\eta=\eta(s)>0$ tel que l'une des deux in\'egalit\'es

\item[$(i)$]
$
\varepsilon^s \le {\varepsilon'}^s N \left( B(x,\varepsilon), \varepsilon' \right),
$
\item[$(ii)$]
$
\varepsilon^s \ge {\varepsilon'}^s N \left( B(x,\varepsilon), \varepsilon' \right),
$\\
est v\'erifi\'ee pour tout $x \in E $
 et pour tout $\varepsilon' \le \varepsilon \le \eta,\, \varepsilon,\varepsilon' \in val(d).$
\end{itemize}
Alors,
$$
\dim_H E= \dimi E.
$$
\end{lemme}

\begin{rem}
Si $(E,d)$ est un espace ultra-m\'etrique compact on peut montrer facilement que
$val(d)$ est l'ensemble des valeurs d'une suite
 $(\varepsilon_n)_{n \ge 0}$
strictement d\'ecroissante et tendant vers $0$. Les conditions $(i)$ et $(ii)$ du lemme pr\'ec\'edent
peuvent alors se reformuler de la mani\`ere suivante~:
Pour tout $s>0$ il existe un entier $K=K(s)$ tel que l'une des deux in\'egalit\'es
\begin{itemize}
\item[$(i)$]
$N(B(x,\varepsilon_n),\varepsilon_{n+1}) \ge \left(\frac{\varepsilon_n}{\varepsilon_{n+1}}\right)^s$
\item[$(ii)$]
$N(B(x,\varepsilon_n),\varepsilon_{n+1}) \le \left(\frac{\varepsilon_n}{\varepsilon_{n+1}}\right)^s$
\end{itemize}
est v\'erifi\'ee pour tout $x \in E$ et pour tout $n \ge K.$
\end{rem}

\begin{prop} \label{propdimHPG}
Si la suite $(\sigma(n))_{n \ge 0}$ est stationnaire ou si elle tend vers l'infini, alors
$$
\dim_H \mathcal{P}(G)=\dimi \mathcal{P}(G).
$$
En particulier, si $G$ est de croissance exponentielle alors
$\dim_H \mathcal{P}(G)=\infty.$
\end{prop}

\begin{cor} \label{propdimHPLm}
Si $m=1$ alors $\dim_H \mathcal{P}(\mathbb{L}_m)= \dim_H \mathcal{P}(\Z)=2 \log 2$.
Si $m \ge 2$ alors
$\dim_H \Pm =\infty.$
\end{cor}

\begin{sdemo}{D\'emonstration de la proposition \ref{propdimHPG}}
On montre d'abord que $\mathcal{P}(G)$ remplit la condition
$(1)$ du lemme \ref{lemmecritere}.
Soit $B$ une boule ferm\'ee de $\mathcal{P}(G)$ de rayon $r=2^{-n} \in val(d_{\mathcal{P}})$ et dont le centre est
 une partie $x$ de $G$.
Si  $B_{X}(n+1)  \nsubseteq x$, on pose $x'=x \cup \{g\}$, o\`u $g$ est un \'el\'ement
 quelconque de $B_{X}(n+1) \backslash x.$
Si $B_{X}(n+1) \subset x$, on pose $x'=B_{X}(n).$
Les parties $x$ et $x'$ sont deux points de $B$ v\'erifiant $d_{\mathcal{P}}(x,x')=r$ si bien que $diam(B)= r$.
Ce qui montre que $\mathcal{P}(G)$ v\'erifie $(1)$.\\
On fixe maintenant $s>0$ et l'on prouve que, si $\sigma$ est croissante ou tend vers l'infini, alors
la condition $(2)$ du lemme \ref{lemmecritere} est aussi remplie.
Soit $n$ un entier positif ou nul.
Alors  $N \left( B(x,2^{-n}),2^{-n-1} \right)
=2^{ \vert B_{X}(n+1) \backslash B_{X}(n) \vert}=2^{\sigma(n)}.$
Si $\sigma$ tend vers l'infini lorsque $n$ tend vers l'infini, il existe $K=K(s)$ tel que
$2^{\sigma(n)} \ge \left( \frac{2^{-n}}{2^{-n-1}} \right)^s=2^s$ pour tout $n \ge K$. Dans ce cas $(2.i)$ est v\'erifi\'ee.
Supposons que $\sigma$ est stationnaire et prend la valeur constante $S<\infty$ pour tout $n \ge K$, pour un certain entier
 $K$. Si $s< \log 2^S$ (respectivement $s \ge \log_2 2^S$) alors $2^{\sigma(n)} \ge 2^s$ (respectivement
 $2^{\sigma(n)} \le 2^S$)
 pour tout $n \ge K$. La condition $(2)$ est donc satisfaite.
Si $G$ est de croissance exponentielle alors $\underset{n \to \infty}{\lim} \sigma(n)=\infty$ (voir \cite[Ch.VI.C, Remarque 53.v]{DeHa})
 et $ \underset{n \to \infty}{\liminf} \frac{\beta(n)}{n}=\infty$. L'application des proposition
 \ref{propdimMPG} et lemme \ref{lemmecritere} conduit \`a l'\'egalit\'e $\dim_H \mathcal{P}(G)=\infty.$
\end{sdemo}
Terminons ce chapitre avec le calcul de la dimension de Minkowski de $\mathcal{G}_1$. Ce calcul est une
application du lemme \ref{lemmerecouvrement2}.
\begin{prop} \label{dimZ}
Lorsque $m=1$, $\mathcal{G}_m=\mathcal{G}(\Z)$ et l'on a
$$
\dims \mathcal{G}(\Z)=0.
$$
\end{prop}

\begin{demo}
Fixons $n \ge 1$ et posons $\varepsilon=2^{-n}$. Par le lemme \ref{lemmerecouvrement2},
$N(\mathcal{G}(\Z),\varepsilon)$ est le cardinal de l'ensemble $\{P \subset \Z \, \vert \, P=k\Z
\cap \{-n,\ldots,-1,0,1,\ldots,n \}
k \in \Z\}$.
Ainsi $N(\mathcal{G}(\Z),\varepsilon)=n+1$, si bien que
$$ \dims \mathcal{G}(\Z)=
\underset{n \to \infty}{\limsup} \frac{\log N(\mathcal{G}(\Z),2^{-n})}{\log 2^n}
=\underset{n \to \infty}{\limsup} \frac{\log (n+1)}{n}=0.$$
\end{demo}
Cette derni\`ere proposition est aussi un cas particulier du th\'eor\`eme \ref{thmdimZm} d\'emontr\'e au chapitre
\ref{secdimGZ}.
\section{Estimations des dimensions de Minkowski d'un sous-espace de groupes \`a petite simplification}\label{secdimPS}

On suppose $m \ge 2$. Un mot r\'eduit $w=avb \in \Lm$ avec $a,b \in S \cup S^{-1}$ est cycliquement r\'eduit si
$a \neq b^{-1}$. Pour tout entier $n \ge 1$, $cyc(n)$ d\'esigne l'ensemble des \'el\'ements de $\Lm$ qui sont
cycliquement r\'eduits et de longueur $n$. On trouve facilement

\begin{equation} \label{eqbn}
\vert B_S(n) \vert=\frac{m}{m-1}((2m-1)^n-1)+1 \, (n \ge 1).
\end{equation}

\begin{equation} \label{eqasympcyc}
\lim_{n \to \infty}\frac{\vert cyc(n) \vert}{(2m-1)^n}=1.
\end{equation}

Soit $u$ et $v$ deux mots cycliquement r\'eduits.
Un sous-mot d'un conjugu\'e cyclique de $u$ ou de $u^{-1}$ qui est aussi un sous-mot de l'un des conjugu\'es cycliques de
 $v$ ou de $v^{-1}$ est appel\'e une \it{pi\`ece} \rm entre $u$ et $v$. Soit $\lambda>0$. Si $R$ est une partie de $\Lm$ form\'ee de mots
 cycliquement r\'eduits telle que toute pi\`ece $p$ entre deux \'el\'ements de $R$ v\'erifie
 $\vert p \vert_S< \lambda \underset{r \in R}{\min} \vert r \vert_S$, on dit que $R$ v\'erifie \it{la condition de petite simplification
 $C'(\lambda)$}. \rm

On suppose dor\'enavant que $\lambda \in \rbrack 0, 1/6 \rbrack$. Pour tout couple d'entiers $n$ et $k$,
on \'ecrit $ps(n)=ps_{k,\lambda}(n)$ l'ensemble
des $k$-uplets d'\'el\'ements de $cyc(n)$
v\'erifiant l'hypoth\`ese m\'etrique de petite simplification $C'(\lambda)$. On d\'efinit les sous-espaces $$PS(n)=\{\cln{r_1,\ldots,r_k} \lhd \Lm \, \vert \, (r_1,\ldots,r_k) \in ps(n)\}
\mbox{ et } PS= \bigcup_{n \ge 1} PS(n).$$

\begin{thm} \label{thmdimPS}
$$
k \log_2(2m-1) \le \dimi PS \le \dims PS \le \frac{k}{1-3\lambda}\log_2(2m-1).
$$
\end{thm}

Le sous-espace $PS$ est d\'enombrable et peut \^etre regard\'e comme l'ensemble
des valeurs d'une suite dans $\Gm$ dont la seule valeur d'adh\'erence
est le groupe distingu\'e trivial $\{1\}$.
 La majoration de la dimension sup\'erieure de Minkowski de cet ensemble
  est rendue possible gr\^ace \`a un th\'eor\`eme classique de la th\'eorie de la petite
 simplification, le th\'eor\`eme de Greendlinger~:

\begin{theor} \cite[Ch.V,Th.4.4]{lyschu} \label{thgreen}
Soit $\mathcal{N}$ la cl\^oture normale d'un ensemble de relateurs $R \subset \Lm$ v\'erifiant la condition $C'(\lambda),
\lambda \in \rbrack 0, 1/6 \rbrack$.
Alors tout \'el\'ement non trivial de $\mathcal{N}$ contient un sous-mot $s$ d'un conjugu\'e cyclique
d'un \'el\'ement $r \in R$ ou de $r^{-1}$
avec $\vert s \vert_{S}>(1-3 \lambda) \vert r \vert_{S}$.
\end{theor}
Le th\'eor\`eme pr\'ec\'edent montre en effet que lorsque les
relateurs de groupes \`a petite simplification sont suffisamment
grands, les cl\^otures normales de ces relateurs ne sont plus
$\varepsilon$-distinguables du sous-groupe
trivial $\{1\}$ pour $\varepsilon>0$ fix\'e. 
Cette consid\'eration permet alors de majorer tr\`es simplement $P(PS,\varepsilon)$.

Minorer la dimension de Minkowski inf\'erieure de $PS$ s'obtient en minorant $P(PS,\varepsilon)$.
Cette minoration repose sur trois observations. Tout d'abord, des cl\^otures normales
distinctes de relateurs de longueurs suffisamment petites sont $\varepsilon$-distinguables. En suite, on sait par un r\'esultat
 de Greendlinger que la cl\^oture normale de relateurs v\'erifiant $C'(\lambda)$ d\'etermine presque uniquement le choix
de ces relateurs. Enfin, un r\'esultat classique de g\'en\'ericit\'e affirme que des $k$-uplets de tels relateurs ayant tous
 m\^eme longueur sont asymptotiquement aussi nombreux que les $k$-uplets de mots cycliquement r\'eduits ayant tous m\^eme
longueur.

\begin{sdemo}{D\'emonstration du th\'eor\`eme \ref{thmdimPS}}
D\'emontrons en premier lieu que
$$
\dims PS \le \frac{k}{1-3\lambda}\log_2(2m-1).
$$
Pour $n \in \N$, fixons $\varepsilon=2^{-n}$.
Remarquons alors que les points de $\bigcup_{j>\frac{n}{1-3 \lambda}} PS(j)$
ne sont pas $\varepsilon$-distinguables.
En effet, soit $\cln{r_1,\ldots,r_k} \in \bigcup_{j> \frac{n}{1-3\lambda}}PS(j)$,
 le th\'eor\`eme \ref{thgreen} implique que
$
\cln{r_1,\ldots,r_k} \cap \B{n}=\{1\}
$
et donc \\ $d_m(\cln{r_1,\ldots,r_k},\{1\}) \le \varepsilon$. Puisque $d_m$ est ultra-m\'etrique, la distance entre deux points
quelconques de $\bigcup_{j>\frac{n}{1-3 \lambda}}PS(j)$ est inf\'erieure ou \'egale \`a
$\varepsilon$.
En conservant les notations du chapitre \ref{secdim}, il s'ensuit que
$$
P(PS,\varepsilon) \le \vert \B{ \frac{n}{1-3 \lambda}} \vert^k \le
\left( \frac{m}{m-1} \right)^k (2m-1)^{\frac{kn}{1-3 \lambda}}.
$$
D'o\`u~:
$
\underset{n \to \infty}{\limsup} \frac{ \log P(PS,2^{-n})}{n \log 2} \le \frac{k}{1-3 \lambda} \log_2(2m-1).
$

Prouvons maintenant la minoration.
 Lorsque $\varepsilon=2^{-n}$, les points de $PS(n)$ sont tous
$\varepsilon$-distinguables. En effet, si $\cln{r_1,\ldots,r_k}$ et $\cln{r'_1,\ldots,r'_k}$ sont deux \'el\'ements
de $PS(n)$ et si l'on suppose en outre que \\
$
d_m(\cln{r_1,\ldots,r_k},\cln{r'_1,\ldots,r'_k}) \le \varepsilon,
$
il s'ensuit que $$\cln{r_1,\ldots,r_k}\cap \B{n} =\cln{r'_1,\ldots,r'_k} \cap \B{n},$$ si bien que
$
\cln{r_1,\ldots,r_k}=\cln{r'_1,\ldots,r'_k}.
$
Majorons maintenant le cardinal d'une fibre de l'application
$$
\begin{array}{ccc}
ps(n)& \longrightarrow & PS(n)\\
(r_1,\ldots,r_k)   & \longmapsto      & \cln{r_1,\ldots,r_k}
\end{array}
$$
  Il est prouv\'e dans \cite{gre} que si l'on a
  $$\cln{r_1,\ldots,r_k}=\cln{t_1,\ldots,t_k} \mbox{ avec }
 (r_1,\ldots,r_k),(t_1,\ldots,t_k) \in ps(n)$$
 alors pour tout $i=1,\ldots,k$, $r_i$ est un conjugu\'e cyclique de l'un des $t_j$ ou de son inverse.
Donc le cardinal d'une fibre est major\'e par $k!(2n)^k$, de sorte que
$
\vert PS(n) \vert \ge \frac{\vert ps(n) \vert }{k!(2n)^k}.
$
 On observe ainsi que
$P(PS,\varepsilon) \ge \frac{\vert ps(n) \vert }{k!(2n)^k}.$
Puisque
\begin{equation} \label{asymsc}
\underset{n \to \infty}{\lim} \frac{\vert ps(n) \vert}{(2m-1)^{kn}}=1,
\end{equation}
 par le lemme 10 de \cite{arzhan}(voir aussi
 \cite{champ1}),
on en d\'eduit que
$$
\underset{\varepsilon \to 0}{\liminf} \frac{ \log P(PS,\varepsilon)}{\log 1/ \varepsilon}=
\underset{n \to \infty}{\liminf} \frac{\log \frac{ps(n)}{k!(2n)^k}}{n \log 2}
 \ge \underset{n \to \infty}{\liminf} \frac{\log \frac{(2m-1)^{kn}}{k!(2k)^n}}{n \log 2}=k \log_2(2m-1).
$$

\end{sdemo}

\begin{cor} Si $m \ge 2$ alors
$\dimi \Gm= \infty.$
\end{cor}

\section{Estimations des dimensions de Minkowski d'un sous-espace de groupes  \`a un relateur } \label{secdimPS}%

On suppose $m \ge 2$. Soit $q \ge 2$ un entier.
On d\'efinit $UR(n)= \{\cln{r^q},\, r \in cyc(n) \}$ et $UR=\bigcup_{n \ge 1} UR(n).$
A nouveau, le sous-espace d\'enombrable $UR$ est l'ensemble des valeurs d'une suite dans $\Gm$ dont
la seule valeur d'adh\'erence est $\{1\}$.

\begin{thm} \label{thmdimUR}
$\frac{\log_2(2m-1)}{q} \le \dimi UR \le \dims UR \le \frac{\log_2(2m-1)}{q-1}.$
\end{thm}

Le principe de majoration est le m\^eme que celui employ\'e dans la preuve du th\'eor\`eme~\ref{thmdimPS}.
Il s'agit d'une majoration du nombre d'\'el\'ements de $UR$ qui sont $\varepsilon$-distinguables. Cette majoration repose
sur un th\'eor\`eme de Neumann, analogue pour les groupes \`a un relateur avec torsion du th\'eor\`eme de Greendlinger.

\begin{theor} \cite[Chap.II, Pr.5.28]{lyschu} \label{thneum}
 Soit un entier $q \ge 2$ et soit
$r$ un mot cycliquement r\'eduit de $\Lm$. On consid\`ere le groupe \`a
un relateur $G= \left< S \vert \, r^q \right>$. Si $w$ est un mot r\'eduit de $\Lm$ dont l'image est triviale
dans $G$, alors il existe un sous-mot r\'eduit $u$ de $w$
qui est aussi un sous-mot de $r^q$ ou de son inverse et v\'erifiant
$\vert u \vert_{S}>(q-1) \vert r \vert_{S}$.
\end{theor}

La minoration des dimensions provient d'une minoration du nombre d'\'el\'ements $\varepsilon$-distinguables et
tient en ces deux observations~: des cl\^otures normales distinctes de relateurs suffisamment
petits sont $\varepsilon$-distinguables; deux mots cyliquement r\'eduits et tels qu'eux m\^eme et leurs inverses
ne sont pas conjugu\'es ont des cl\^otures normales distinctes. Ce dernier fait est un r\'esultat de Magnus.

\begin{sdemo}{D\'emonstration du th\'eor\`eme \ref{thmdimUR}}
Pour $ n \in \N$, fixons $\varepsilon=2^{-n}$.
Montrons que les points de $\bigcup_{i> \frac{n}{q-1}} UR(i)$ ne sont pas $\varepsilon$-
distinguables. Si $\cln{r^q}$ est tel que $r \in cyc(i)$ avec $i > \frac{n}{q-1}$, alors
$d_m(\cln{r^q}, \{1\}) \le \varepsilon$. En effet, le th\'eor\`eme
\ref{thneum} assure que tout mot
$w \in \cln{r^q}$ est tel que $\vert w \vert_{S}>(q-1) \vert r \vert_{S}$, si bien que
$ \cln{r^q} \cap \B{n}=\{1\}$. Puisque $d_m$ est ultra-m\'etrique, la distance entre deux points
quelconques de $\bigcup_{i > \frac{n}{q-1}} UR(i)$ est inf\'erieure ou \'egale \`a
$\varepsilon$. Ainsi
 \begin{equation} \label{ineqsup}
 P(UR,\varepsilon) \le \vert \bigcup_{i \le \frac{n}{q-1}}UR(i) \vert
  \le \vert \B{\frac{n}{q-1}} \vert \le \frac{m}{m-1}(2m-1)^{\frac{n}{q-1}}.
 \end{equation}
Montrons maintenant que les points de $\bigcup_{i \le \frac{n}{q}}UR(i)$ sont $\varepsilon$-distinguables.
En effet, si $\cln{r^q}$ et $\cln{r'^q}$ sont des \'el\'ements de $\bigcup_{i \le \frac{n}{q}}UR(i)$
tels que \\ $d_m(\cln{r^q},\cln{r'^q}) \le \varepsilon $, alors \\
$\cln{r^q} \cap \B{n} =\cln{r'^q} \cap \B{n}$, si bien que
$\cln{r^q}=\cln{r'^q}$ puisque $\vert r^q \vert_{S}, \vert r'^q \vert _{S}\le n.$

Par un r\'esultat d\^u \`a Magnus \cite[Chap.II, Pr. 5.8]{lyschu}, on sait que lorsque
$r_1$ et $r_2$ sont des \'el\'ements cycliquement r\'eduits de $\Lm$ tels que
 $\cln{r_1}=\cln{r_2}$, alors $r_1$ est un conjugu\'e cyclique de $r_2$ ou de $r_2^{-1}.$
Le cardinal d'une fibre de l'application
$$
\begin{array}{ccc}
cyc( \lbrack \frac{n}{q} \rbrack) & \rightarrow & UR( \lbrack \frac{n}{q} \rbrack) \\
  r   & \mapsto     & \cln{r^q}
\end{array}
$$
est donc au plus $2 \lbrack \frac{n}{q} \rbrack $, si bien que $\vert \bigcup_{i \le \frac{n}{q}}UR(i) \vert
\ge \vert UR(\frac{n}{q})\vert
\ge \frac{\vert cyc( \lbrack \frac{n}{q} \rbrack ) \vert}{2 \lbrack \frac{n}{q} \rbrack}$.
 Ainsi,
\begin{equation} \label{ineqinf}
  P(UR,\varepsilon) \ge \frac{\vert cyc( \lbrack \frac{n}{q} \rbrack) \vert}{2 \lbrack \frac{n}{q} \rbrack}.
\end{equation}
Le th\'eor\`eme s'en d\'eduit par passage \`a la limite lorsque
$n$ tend vers l'infini dans les in\'egalit\'es (\ref{ineqsup}) et (\ref{ineqinf}) o\`u l'on applique (\ref{eqasympcyc}).
\end{sdemo}

\section{ Dimension de Minkowski sup\'erieure du Cantor
$\mathfrak{B}$ de Grigorchuk} \label{secdimB}

Dans ce chapitre on montre que la dimension de Minkowski sup\'erieure de l'espace
de Cantor construit par Grigorchuk dans \cite[Ch.6]{gri} est nulle. Rappelons cette construction.
 Grigorchuk d\'efinit \`a partir d'un algorithme une famille
 de groupes $\widetilde{G}_{\omega}$ param\'etr\'es par des suites $\omega$ de symboles
 pris dans l'ensemble $\{\overline{0},\overline{1},\overline{2}\}$.
  Une sous-famille de ces groupes est form\'ee de groupes de croissance interm\'ediaire qui r\'epondent \`a une question
  pos\'ee par Milnor. Grigorchuk montre aussi \cite[Ch.6, Pr. 6.2]{gri}
  que l'ensemble de tous ces groupes est un espace de Cantor lorsqu'on le munit de la topologie de Cayley.\vskip 0.5 cm
 On fixe l'alphabet $\{a,b,c,d\}$ et l'on d\'efinit $\ov{0}=(a,a,1), \ov{1}=(a,1,a)$ et $\ov{2}=(1,a,a)$ o\`u $1$ d\'esigne
 le mot vide.
On d\'esigne par $\mathbb{F}(a,b,c,d)$ le groupe libre sur $\{a,b,c,d \}$, par $\vert w \vert$ la longueur d'un mot
$w \in \mathbb{F}(a,b,c,d)$ relativement \`a $\{a,b,c,d\}$ et par
$\Omega$ l'ensemble des suites d'\'el\'ements de $\{\ov{0},\ov{1},\ov{2} \}$.

 Pour tout mot $w \in \mathbb{F}(a,b,c,d)$, on consid\`ere la forme r\'eduite positive $r(w) \in \mathbb{F}(a,b,c,d)$
  obtenue en appliquant les r\`egles de r\'e\'ecriture suivantes :
 \begin{itemize}
 \item[$(i)$] $x^{-1} \rightarrow x$
 \item[$(ii)$] $x^2 \rightarrow 1$
 \item[$(ii)$] $xy \rightarrow z$
 \end{itemize}
 o\`u $x,y,z \in \{a,b,c,d\}$ pour les deux premi\`eres r\`egles et $x,y,z \in \{b,c,d\}$ et sont distincts
  pour la troisi\`eme. L' applications de ces r\`egles \`a un mot $w$ sont appel\'ees \it
  simplifications \'el\'ementaires. \rm It\'er\'ees jusqu'au moment o\`u plus aucune r\`egle ne s'applique,
  ces simplifications donnent une forme r\'eduite positive $r(w)$ de $w$ dans le groupe
  $$
  \Gamma= \left< a,b,c,d \,\vert \, a^2=b^2=c^2=d^2=bcd=1 \right> \simeq \Z_2 \ast (\Z_2 \times \Z_2).
  $$
 Ainsi, $w \overset{\Gamma}{=}1$ si et seulement si $r(w)=1$.

 Ayant fix\'e $\omega \in \Omega,\, \omega=(a_n,b_n,c_n)_{n \ge 1}$, on d\'efinit comme dans \cite[Ch.6, pp 287-288]{gri},
 l'ensemble $\widetilde{S}_{\omega}$ des \'el\'ements de $\mathbb{F}(a,b,c,d)$ pour lesquels
 l'algorithme $\mathfrak{a}$ d'oracle $\omega$  (\cite[Ch.2]{gri} et \cite[pp 84-86]{gri2})
 que l'on d\'ecrit plus bas, conduit \`a un r\'esultat positif.
 Pour tout mot positif $w=w(a,b,c,d)$  o\`u la lettre $a$ apparait un nombre pair de fois, on
  d\'efinit deux proc\'ed\'es de r\'e\'ecriture $\varphi_0^{(n)}$ et $\varphi_1^{(n)}$ pour tout $n \in \N^{\ast}$.
  Le r\'esultat $\varphi_i^{(n)}(w)$ du $i$-\`eme processus de r\'e\'ecriture, $i=0,1$, est un mot sur l'alphabet
  $\{a,b,c,d\}$ obtenu \`a partir de $w$ en associant \`a chaque lettre de $w$ un symbole de l'alphabet
  $\{a,b,c,d \}$ en observant les r\`egles de substitutions suivantes~:
  $$
  \varphi_i^{(n)} : \left\{
  \begin{array}{ccc}
  a & \rightarrow & 1 \\
  b & \rightarrow & b_n \\
  c & \rightarrow & c_n \\
  d & \rightarrow & d_n
  \end{array}
  \right.
  $$
  si le nombre de lettres $a$ dans $w$ pr\'ec\'edant le symbole courant auquel on applique la r\`egle de substitution
  est pair lorsque $i=0$ ou impair lorsque $i=1$. De mani\`ere analogue,
  $$
  \varphi_i^{(n)} : \left\{
  \begin{array}{ccc}
  a & \rightarrow & 1 \\
  b & \rightarrow & b \\
  c & \rightarrow & c \\
  d & \rightarrow & d
  \end{array}
  \right.
  $$
  si le nombre de lettres $a$ dans $w$ pr\'ec\'edant le symbole courant auquel on applique la r\`egle de substitution
  est impair lorsque $i=0$ ou pair lorsque $i=1$.
  L'algorithme $\mathfrak{a}$ d'oracle $\omega$ est le suivant~:
  \begin{algo} Pour d\'ecider si $\widetilde{S}_{\omega}$ contient le mot $w=w(a,b,c,d)$ on proc\`ede comme suit.
  \begin{itemize}
  \item[$(1)$]
  On \'evalue la somme des exposants de la lettre $a$ dans $w$. Si cette somme est impaire alors $w$ n'est pas retenu.
  Si elle est paire on calcule $r(w)$. Si $r(w)=1$, alors $w$ est retenu et l'algorithme s'arr\`ete.
   Si $r(w) \ne 1$
  et si $ \vert r(w) \vert=1$ alors $w$ n'est pas retenu et l'algorithme s'arr\`ete. Si $\vert r(w) \vert>1$,
  on proc\`ede \`a l'\'etape $(2)$
  \item[$(2)$]
  On calcule $w_i=\varphi_i^{(1)}(r(w)),i=0,1,$ et l'on retourne \`a l'\'etape $(1)$
  o\`u les v\'erifications s'appliquent
  maintenant \`a deux mots $w_i,i=0,1.$ Si l'algorithme se poursuit apr\`es $2n$ \'etapes, dans l'\'etape $(1)$ les
  v\'erifications portent alors sur $2^n$
  mots $$w_{0 \ldots 0},\dots,w_{i_1 i_2 \ldots i_n}, \ldots ,w_{1 \ldots 1}$$ et l'on calcule
   ensuite dans l'\'etape $(2)$, si elle a lieu,
    $$\varphi_i^{(n+1)}(r(w_{0 \ldots 0})),
   \dots,\varphi_i^{(n+1)}(r(w_{i_1 i_2 \ldots i_n})), \ldots ,\varphi_i^{(n+1)}( r(w_{1 \ldots 1})),\quad i=0,1.$$
  \end{itemize}
  \end{algo}

  Puisque $\vert \varphi_i^n(r(w)) \vert \le \frac{\vert w \vert+1}{2}$ pour $i=0,1$ et tout $n \ge 1$, un mot $w$ est
  accept\'e ou rejet\'e au bout d'au plus $\log_2(\vert w \vert)$ applications r\'ecursives de cet algorithme.
  L'ensemble $\widetilde{S}_{\omega}$ est l'ensemble des mots $w=w(a,b,c,d)$ retenus par l'algorithme.
   On v\'erifie ais\'ement que $\widetilde{S}_{\omega}$ est un sous-groupe distingu\'e de $\mathbb{F}(a,b,c,d)$ contenant
   les relateurs qui d\'efinissent $\Gamma$.

  Le groupe marqu\'e $\widetilde{G}_{\omega}$ est le quotient $\mathbb{F}(a,b,c,d)/\widetilde{S}_{\omega}$ marqu\'e par l'image
  des g\'en\'erateurs $\{a,b,c,d\}$. Pour simplifier on identifie $\mathcal{G}_4$ et $\mathcal{G}(a,b,c,d)$.

Un algorithme compl\'ementaire d\'ecrit dans \cite[Ch. 5]{gri},
 permet de reconstituer les $n$ premiers symboles de $\omega$ en connaissant les \'el\'ements
 de $\widetilde{S}_{\omega}$ de longueur inf\'erieure ou \'egale \`a $2^{n+2}$.
Les deux algorithmes consid\'er\'es par Grigorchuk conduisent aux r\'esultats de la proposition 6.2
dont une reformulation est la suivante :

\begin{prop} \rm{ \cite[Ch.6]{gri}} \label{propcantor}
\begin{itemize}
\item[$(i)$]
Soit $\omega_1, \omega_2 \in \Omega$ et soit $n \in \N$.
Si $\omega_1$ et $\omega_2$ coincident sur les $n$ premi\`eres coordonn\'ees, alors
$d_4(\widetilde{S}_{\omega_1},\widetilde{S}_{\omega_2}) \le 2^{-2^n}$.
\item[$(ii)$]
Si $d_4(\widetilde{S}_{\omega_1},\widetilde{S}_{\omega_2}) \le 2^{-2^{n+2}}$, alors
$\omega_1$ et $\omega_2$ coincident sur les $n$ premi\`eres coordonn\'ees.
\item[$(iii)$]
Le sous-espace $\mathfrak{B}=\{\widetilde{S}_{\omega} \}_{ \omega \in \Omega}$ de
$\mathcal{G}_4$ est un espace de
Cantor.
\end{itemize}
\end{prop}

\begin{thm} \label{thmdimBeta}
La dimension de Minkowski sup\'erieure du sous-espace $\mathfrak{B} \subset \mathcal{G}_4$
est nulle.
\end{thm}

\begin{sdemo}{D\'emonstration du th\'eor\`eme \ref{thmdimBeta}}
Pour $n \in \N$, fixons $\varepsilon=2^{-2^n}$.
 On consid\`ere dans $\mathcal{G}_4$ l'ensemble des boules de rayon $\varepsilon$ centr\'ees en les groupes
$\widetilde{S}_{\omega}$ pour lesquels $\omega$ est constante \`a partir de la $n$-\`eme coordonn\'ee. Il s'agit
d'un recouvrement de $\mathfrak{B}$ en vertu du point $(i)$ de la proposition \ref{propcantor}. Le cardinal de ce
recouvrement est au plus $3^{n}$. Ainsi,
$\frac{\log N(\mathfrak{B}, \varepsilon)}{\log 1/ \varepsilon} \le \frac{\log 3^n}{2^n \log 2}$.
D'o\`u $\dims \mathfrak{B}=0$.
\end{sdemo}

\section{Dimension de Minkowski sup\'erieure du sous-espace des groupes commutatifs marqu\'es
 \`a $m$ g\'en\'erateurs} \label{secdimGZ}

Soit $\mathcal{A}=(a_1,\ldots,a_m )$ l'ensemble ordonn\'e des
 g\'en\'erateurs canoniques du groupe ab\'elien
libre $\Z^m$. On d\'esigne  par $\BZn$ l'ensemble des \'el\'ements de $\Z^m$ de longueur inf\'erieure
ou \'egale \`a $n$ relativement \`a $\mathcal{A}$.
Soit $b_n$ le cardinal de $\BZn$. Il est imm\'ediat de v\'erifier que $\underset{n \to \infty}{\lim} \frac{b_n}{n^m}
=vol(K)$, o\`u $vol(K)=\frac{2^m}{m!}$ est le volume de l'enveloppe convexe $K$ de $\mathcal{A} \cup \mathcal{A}^{-1}$ dans $\R^m$.
Il existe donc $C \ge 0$ tel que $b_n \le C n^m$ pour tout $n \in \N$.

 On conserve les notations de l'introduction en consid\'erant l'espace $\mathcal{G}(\Z^m)$
 muni de la m\'etrique $d_{\mathcal{G}}$. L'\'epimorphisme de $\Lm$ sur $\Z^m$ d\'efini \`a partir de la bijection naturelle entre les ensembles de g\'en\'erateurs
 ordonn\'es $S$ et $\mathcal{A}$ induit un plongement de $\mathcal{G}(\Z^m)$ dans $\mathcal{G}_m$ qui
 est isom\'etrique sur son image. Cette image est l'ensemble des groupes commutatifs marqu\'es.

\begin{thm} \label{thmdimZm}
La dimension de Minkowski sup\'erieure de l'espace m\'etrique $(\mathcal{G}(\Z^m),d_{\mathcal{G}})$ est nulle.
\end{thm}

\begin{demo}
Soit $n$ un entier sup\'erieur ou \'egale \`a deux. Soit $\mathcal{P}_n=\{ \BZn \cap R \,\vert \, R \le \Z^m \}$.
Par le lemme \ref{lemmerecouvrement2}, $N(\GZm,2^{-n})= \vert \mathcal{P}_n \vert$.
Le cardinal de $\mathcal{P}_n$ est inf\'erieur au nombre de sous-groupes de $\Z^m$ dont une base est dans la boule
$B_{\mathcal{A}}(n)$.
On a donc $\vert \mathcal{P}_n \vert \le \sum_{l=0}^m \binom{b_n}{l} \le b_n^{m+1}.$
Il s'ensuit que
 $$\dims \GZm =\limsup_{n \to \infty }\frac{\log N(\GZm,2^{-n})}{n} \le \limsup_{n \to \infty}
  \frac{\log n^{m(m+1)}}{n}=0.$$
\end{demo}

\end{document}